

\input vanilla.sty
\input definiti.tex
\magnification 1200
\baselineskip 18pt
\overfullrule=0pt
\input mathchar.tex

\def\lam{\lambda}

\def\gam{\gamma}

\def\Om{\Omega}
\def\om{\omega}

\def\varp{\varphi}

\def\vare{\varepsilon}

\def\Sig{\Sigma}
\def\suml{\sum\limits}

\def\part{\partial}

\def\Vol{\text{\rm Vol}}

\def\tilg{\tilde g}

\par\midspace{2cm}\par

\title
YAMABE SPECTRA
\endtitle
\author
Alexander G. Reznikov\footnote"$^{*}$"{Partially sponsored by
the Edmund Landau Center for research in Mathematical Analysis, supported by
the Minerva Foundation (Germany).}\\
Institute of Mathematics \\
The Hebrew University \\
Givat Ram 91904
ISRAEL
\endauthor

\centerline{Preprint No. 16}
\centerline{1992/93}
\vskip40pt

\vfill\eject

\heading{Contents}\endheading

{\settabs 5\columns
\+0. Introduction & \dotfill&\dotfill&\dotfill&3\cr
\+1. The Homological Plateau Problem & &\dotfill&\dotfill&5\cr
\+2. The Thurston norm and the area norm & &&\dotfill&6\cr
\+3. Geometry of Jacobians and conformal invariants
& &&\dotfill&8\cr
\+4. Obstructions to Nash embeddings & &&\dotfill&10\cr}

\vfill\eject

\title
Yamabe Spectra
\endtitle
\author
Alexander G. Reznikov
\endauthor

\heading{0. Introduction}\endheading

The Uniformization Program of William Thurston prescribes, for a large class of
compact closed three-manifolds, the existence of a hyperbolic metric, i.e. a
metric
of constant negative
curvature - $1$. One approach to find such a metric
is a two-step variational problem, which we now describe.
Consider a compact Riemannian manifold $(M^n,g)$, and for any metric $\tilde g$
in the conformal class $C=C(g)$, define
$$
Y(\tilde g)=
\frac{\int_M R_{\tilde g}d\Vol_{\tilde g}}{[\Vol(\tilg)]^{\frac{n-2}{n}}}
\tag 1
$$
Then one easily shows that the  critical points $\tilde g$ of $Y(\tilg)$ in
$C$ are metrics of constant scalar curvature [16]. The minimum
$\inf Y(\tilg)$, $\tilg\in C$ is called the Yamabe invariant of $(M,g)$ and
denoted
$\lam(g)$. The solution to the Yamabe problem,
obtained in the  last thirty years due to the efforts of Yamabe [16], Trudinger
[14],
Aubin [1], and Schoen [11], [12] establishes the existence of the
minimizing metric $\tilg$ in $C$.
For $M$ three-dimensional one knows, under certain topological restrictions,
that
the constant scalar curvature of a metric on $M$ may not be positive.
In fact, hyperbolizable manifolds never carry a metric of nonnegative scalar
curvature, as shown by Schoen-Yau [13] and Gromov-Lawson [3].
So for such manifolds, $\lam(g)$ is always negative.
The importance of the functional
$C\to\lam(g)$ stems from the hope to find a conformal class $C$, maximizing
this functional.
Then an easy computation show that the corresponding constant scalar curvature
metric
$\tilg$ in $C$ is actually Einstein, hence hyperbolic, since $\dim M=3$.
It is more convenient to introduce a modified function
$Vol:C\mapsto Vol(C)$, as follows:
take the appropriately scaled metric $\tilg$ in $C$ of constant scalar
curvature $-1$,
and denote $Vol(C)=\Vol(\tilg)$.
We will call the set
$\{Vol(C)\}$ the Yamabe spectrum of $M$.

The expectation for the global  minimum of $Vol$ at the hyperbolic metric is
justified
by the following known facts: firstly, if $g_0$ is a (necessarily unique up to
a diffeomorphism)
hyperbolic metric on $M$, then $(D^2Vol)_{g_0}>0$, see [4], theorem 8.2.
Secondly, if $\dim M=4$ and $g_0$ is hyperbolic, then the Euler characteristic
and  signature
computations of Johnson and Millson [4], theorem 8.3 show that $Vol$ attains
its global minimum at $g_0$.

In the three-dimensional case we do not have Gauss-Bonnet type formulas, and
one
should look for other ways for estimation $Vol$.
In the present paper, we deal with Haken three-manifolds with infinite first
homology group.
For every conformal class, $C$, we introduce an invariant of $C$, coming from
the
$L^3$-geometry of the Jacobian variety
$J^1(M)=H^1(M,\bbr)/H^1(M,\bbz)$.
We prove that this invariant estimates $Vol(C)$ from below and use it to show
our
main result:

\proclaim{Theorem 2} Let $M$ be a compact oriented homologically atoroidal
three-manifold with $\pi_2(M)=0$ and $H_1(M,\bbr)\not=0$.
The $\sup\limits_CVol(C)=\infty$.
\endproclaim

We use our estimates for demonstrating  new global obstructions for the Nash
isometrical immersions $M^3\to N$ of arbitrary codimension,
which induce nontrivial map in the first homology.
Finally, we establish, for all three-manifolds with the pinched negative
curvature
$-K\le K(M)\le -k<0$, the following principle:
the Thurston genus norm in $H_2(M)$ is uniformly equivalent to the area norm.

This research was started when I was visiting the Ruhr-Universit\"at Bochum.
I wish to thank J. Jost and S. Wang for very helpful discussions.

\vfill\eject

\heading{1. The Homological Plateau Problem}\endheading

The following version of the existence theorem for minimal surfaces
of Sacks-Uhlenbeck-Schoen-Yau (c.f. [10], [13]) has been established recently
by
Marina Ville [16]:

\proclaim{Theorem A} Let $M$ be a compact Riemannian manifold and let
$0\not=z\in H_2(M,\bbz)$. Then there exists a (branched) minimal
immersion $f$ of a compact oriented, possibly disconnected Riemannian surface
$\Sig^g$ to $M$, representing $z$.
Moreover, one can take $\Sig^g$ of the least possible genus among all surface,
representing $z$ and the map $f$ to be globally area minimizing in its homology
class.
\endproclaim

For the reader's convenience we sketch the proof here.
Start with the following lemma.

\proclaim{Lemma 1} (comp. [8], [17]). Let $X$ be a $CW$-complex and let
$z\in H_2(X,\bbz)/\overline{\pi_2(X)}$. There exists a collection of oriented
Riemannian
surfaces $\Sig^{g_1},\ldots,\Sig^{g_k}$, $g_i\ge1$ and maps $f_i:\Sig^{g_i}\to
X$ such that
\item{\rm (i)} $\sum(2g_i-2)=||z||_g$ (see the  section 2).
\item{\rm (ii)} $\suml_i[f_i]=z$
\item{\rm (iii)} For any essential simple loop $\gam$ in $\Sig^{g_i}$,
$f_{i_*}([\gam_i])\not=0$ in $\pi_1(x)$.
\endproclaim

\demo{Proof} Start with any collection of $\Sig^{g_i}$ satisfying (i) and (ii).
Suppose (iii) is not valid, that is, for some $\gam$ in $\Sig^{g_i}$,
there exists a map $\varp:D^2\to X$ with $\varp|\part D^2=f|\gam$
(we identify $\gam$ with $\part D^2$).
Cut $\Sig^{g_i}$ along $\gam$ and paste two copies of $D^2$ along the new
boundaries.
Denote $\tilde\Sig_i$ the resulting surface and let $\tilde f_i:\tilde\Sig_i\to
X$
be the map, obtained by patching $\varp$ and $f_i$.
Consider the two following cases.
\item{1)} $\gam$ is not separating. Then genus $(\tilde\Sig_i)=g_i-1$

\item{2)} $\gam$ is separating. Then if $\tilde\Sig^{(1)}_i$ and
$\tilde\Sig^{(2)}_i$
are the two components of $\Sig^{g_i}-\gam$ , then
$\tilde\Sig_i=\tilde\Sig_i^{(1)}\sqcup\tilde\Sig_i^{(2)}$ and genus
$(\tilde\Sig^{(1)})+$ genus $(\tilde\Sig^{(2)})=g_i$.

In any case, we have, first, that the new collection of surfaces and maps still
represents $z$, and, second, that $\sum|\chi(\Sig^{g_i})|$ strictly decreases.
This contradicts (i) and thus (iii) is valid.

Returning to the proof of Theorem A, we first note, that by the theorem of
Sacks and
Ulenbeck, there exists a set of minimal spheres in $M$, generating
$\pi_2(M)$ as a $\pi_1(M)$ module, so their images in
$\overline{\pi_2(M)}$ generate $\overline{\pi_2(M)}$.
Hence we may work modulo $\overline{\pi_2(M)}$.
Fix $z\in H_2(M)/\overline{\pi_2(M)}$ and consider a
collection of surfaces and maps $\Sig^{g_i}$, $f_i$ as in the lemma 1, so that
$\sum(2g_i-2)$ is the least possible.
By the Hopf exact sequence, we have $H_2(M)/\overline{\pi_2(M)}=H_2(\pi_1(M))$.
Now, the condition (iii) implies, by the theorem of Schoen-Yau [13] and
Sacks-Ulenbeck
[10], that there exist, for any $i$, a branched minimal immersion
$\psi_i:\Sig^{g_i}\to M$, inducing the same action
on $\pi_1(\Sig^{g_i})$
as $f_i$, up to a conjugation.
In particular, $\Sig[\psi_i]=z$ in $H_2(\pi_1(M))$, as desired.

We refer to [8] for the further refinement of this result and numerous
algebraic applications.

If the dimension of $M$ is three, we may (and will) consider $f$ to be
unbranched,
by the same argument as in [13].

\heading{2. The Thurston Norm and the Area Norm}\endheading

Let $M$ be a smooth Riemannian manifold and let
$z\in H_2(M,\bbz)/\overline{\pi_2(M)}$
where $\overline{\pi_2(M)}$ is the image of the Hurewitz map.
In [15], Thurston introduced a seminorm $||z||_g$, which we will take in a form
$$
||z||_g=\inf\limits_{[\Sig]=z}|\chi(\Sig)|
$$
taken over all singular surfaces $f:\Sig\to M$ of genus $\ge1$, representing
$z$.
One makes $||z||_g$ to a norm by the standard normalization procedure
(see [15]) and shows that this norm is essentially equivalent to the
Gromov's simplicial norm.

On the other hand, given a metric on $M$, one has the usual area norm,
or the mass, of $z$:
$$
||z||_a=\inf\limits_{[\Sig]=z}\text{area}(\Sig).
$$

The celebrated Thurston's inequality relates $||z||_g$ to $||z||_a$
in the case when $M$ is of negative curvature.
The following sharp version was established in [7]: if
$-K\le K(M)\le -k<0$ and $M$ is compact, then
$$
||z||_a\le\frac{2\pi}{k}||z||_g.\tag 2
$$

The main purpose of this section is to establish the following result,
which is in a sense converse to (2):

\proclaim{Theorem 1} Let $M$ be compact three-dimensional homologically
atoroidal Riemannian manifold, and let $R(M)=\sup\limits_M(-R(x))$.
Then for any
$z\in H_2(M,\bbz)/\overline{\pi_2(M)}$, one has
$$
||z||_a\ge\frac{2\pi}{R(M)}||z||_g\tag 3
$$
which is sharp.
\endproclaim

\demo{Proof} We recall that a three-manifold $M$ is called homologically
atoroidal if any map of a torus $T^2$ to $M$ induces zero homomorphism
in the second homology.
All hyperbolizable manifolds are atoroidal.

We begin with finding a minimal map $f:\Sig^g\to M$ representing $z$
which exists by the Theorem A.
We assume $|\chi(\Sig)|$ is minimal possible and $f$ is an immersion,
($\Sig^g$ may be disconnected). Let
$\Sig^{g_i}$, $i=1,\ldots,q$, be the components of $\Sig$.
Then $g_i>1$ since $\chi(\Sig)$ is minimal and $M$ is atoroidal.
Thus the second variation formula gives (comp. [13], ( 5.3))
$$
\suml_i\int_{\Sig^{g_i}}R(x)d\text{area}\le\suml_i\int_{\Sig^{g_i}}Kd\text{area}
$$
which implies by Gauss-Bonnet
$$
R(M)\text{Area}(f)\ge2\pi\sum-\chi(\Sig^{g_i})=2\pi||z||_g.
$$
Since $f$ is the minimizing map, we get $\text{Area}(f)=||z||_a$,
so
$$
R(M)||z||_a\ge2\pi||z||_g,
$$
as prescribed by (3).

\proclaim{Corollary 2} Let $T(k,K)$ denotes a class of compact three-manifolds
with negative curvature satisfying $-K\le K(x)\le -k<0$.
Then the Thurston norm and the area norm are uniformly equivalent in
$T(k,K)$.
More precisely, for any $M\in T(k,K)$ and
$z\in H_2(M,\bbz)$, one has
$$
\frac{1}{k}||z||_g\ge\frac{1}{2\pi}||z||_a\ge\frac{1}{3K}||z||_g\tag 4
$$
and the left hand side is sharp.
\endproclaim

\heading{3. Geometry of Jacobians and Conformal Invariants}\endheading

For a compact Riemannian manifold $N$ we denote
$J^1(N)=H^1(N,\bbr)/H^1(N,\bbz)$.
Let $\Om^*(N)$ stand for the de Rham complex of $N$.
For $p\ge1$ and $\om\in \Om^k(N)$ we denote as usual
$||\om||_{L^p}=(\int_N|\om|_x^pd\text{Vol})^{1/p}$.
This induces a norm on $H^k(N,\bbr)$ by the formula
$||w||_{L^p}=\inf\limits_{w\in w}||\om||_{L^p}$.
We claim:

\proclaim{Proposition 3} Let $M$ be a homologically atoroidal three-manifold
with $\pi_2(M)=0$ and let
$0\not=w\in H_1(M,\bbz)$. Then
$$
||w||_{L^1}\ge\frac{4\pi}{R(M)}.\tag 5
$$
\endproclaim

\proclaim{Corollary 4} The volume of the Jacobian $J^1(M)$ in the $L^2$-metric
is at least\hfill\break
$\left(\frac{4\pi}{R(M)\Vol^{1/2}(M)}\right)^m\Vol B_m$, where
$m=b_1(M)$, and $B_m$ stands for the Euclidean ball.
\endproclaim

\demo{Proof of the Proposition 3} Let $\om\in\Om^1(M)$ with $[\om]=w$.
Since all periods of $\om$ are integers, there exists a smooth map
$\varp:M\to S^1=\bbr/\bbz$ with
$\varp^*(dt)=\om$. Let $S(t)=\varp^{-1}(t)$ and write
$$
\int_M||\nabla \varp||d\Vol =\int_0^1\text{Area}(S(t))dt
$$
by the coarea formula.
For almost all $t$, $S(t)$ is smooth and
$[S(t)]\in H_2(M,\bbz)$ is Poincar\'e dual to $w$.
Applying Theorem 1, we get
$$
||\om||_{L^1}\ge 2\pi R^{-1}(M)||PD(w)||_g\ge\frac{4\pi}{R(M)},
$$
so $||w||_{L^1}=\inf\limits_{[\om]=w}||\om||_{L^1}\ge\frac{4\pi}{R(M)}$.
\hfill Q.E.D.

\demo{Proof of the Corollary 4} This follows readily from (5) and the H\"older
inequality.
Let $M,w$ be as in the proposition 3, and let
$\om\in\Om^1(M)$ with $[\om]=w$.
Let $g$ be the metric of $M$ and let $h$ be its conformal perturbation.
Say $h=\varp\cdot g$ for some positive
$\varp\in C^\infty(M)$.
Using the proposition 3, we get
$$
R_h(M)\cdot\int_M||\om||_h d\Vol_h\ge4\pi,
$$
or
$$
R_h(M)\cdot\int_M||\om||_g\cdot\varp^2d\Vol_g\ge4\pi.
$$
This gives
$$
\int_M\varp^2d\Vol_g\ge\frac{4\pi}{R_h(M)\cdot ||\om||_{L^\infty_g}},
$$
and, by H\"older,
$$
\Vol_h(M)=\int_M\varp^3d\Vol_g\ge
\left(\frac{4\pi}{R_h(M)\cdot||\om||_{L^\infty_g}}\right)^{3/2}\cdot
\Vol_g^{-1/2}(M),
$$
so
$$
R_h^{3/2}(M)\cdot \Vol_h(M)\ge
\left(\frac{4\pi}{||\om||_{L^\infty_g}}\right)^{3/2}
\Vol_g^{-1/2}(M).\tag 6
$$
We wish to improve this, letting the original metric $g$ to change within its
conformal class.
Put $\hat g=\psi\cdot g$ and write (5) for $\hat g$ instead of $g$ to get
$$
(4\pi)^{3/2}R_h^{-3/2}(M)\Vol^{-1}_h(M)\le
(\sup||\om(x)||\cdot\psi^{-1}(x))^{3/2}\left(\int\psi^3d\Vol_g\right)^{1/2}.
$$
The infimum of the right hand side taken over all $\psi>0$ is easily seen
to be $||\om||_{L^3_g}^{3/2}$, so we get finally
$$
4\pi R^{-1}_h(M)\Vol_h^{-2/3}(M)\le||\om||_{L^3_g}.
$$
Letting $h$ be the Yamabe metric in $C_g$, we arrive to the following result

\proclaim{Proposition 5} Let $M$ be compact three-dimensional
homologically atoroidal manifold with $\pi_2(M)=0$ and infinite $H_1(M,\bbz)$.
For any conformal class $C$ we have
$$
Vol(C)\ge\frac{4\pi}{||w||_{L^3}},\tag 7
$$
where $w$ is any class in $H^1(M,\bbz)$ and the $L^3$-norm is taken according
to any metric $g\in C$.
\endproclaim

Observe that $||w||_{L^3_g}$ does not depend on the choice of the metric.
In fact, $L^3$ geometry of the Jacobain $J^1(M)$  depends only on $C$.
The number $\inf\limits_{w\in H^1(M,\bbz)}||w||_{L^3}$ denoted
$j(C)$, is therefore a  conformal invariant.
We can write (6) in the form
$$
Vol(C)\ge 4\pi j^{-1}(C). \tag 8
$$
One can view (8) as a three-dimensional version of the Li-Yau estimates, c.f.
[5].

\proclaim{Theorem 2} Let $M$ be as in the Proposition 5. Then
$$
\sup\limits_CVol(C)=\infty.
$$
\endproclaim

\demo{Proof} Fix $w\in H^1(M,\bbz)$ and $\om\in w$. Fix $\vare>0$.In view of
the proposition 5, it is
enough to find a metric $g$ on $M$ such that $||\om||_{L^3_g}<\vare$.
For that, fix a smooth measure $\mu$ on $M$.
We will always assume that the density $\mu|_g=1$.
Set $||\om||_g<\vare^{1/3}$ everywhere and correct $g_x$ if necessary in the
kernel of $\om_x$ to achieve $\mu|_g=1$, keeping $||\om||_g$  unchanged.
This would be a desired metric.
\proclaim{Remark}
\endproclaim
 We show in [9] that the ``most'' of homology three-spheres are Haken. It would
be very interesting to know if the Theorem 2 is still valid for such manifolds.

\heading{4. Obstructions to Nash Embeddings}\endheading

 Let $M^m$ and $Q^q$ be Riemannian manifolds with $q\ge \frac{m(m+1)}{2}$.
Then any distance-decreasing map $\varp:M\to Q$ can be $C^0$-approximated by
an isometrical embedding (c.f. Gromov [2] for the contemporary survey
of related results). In particular, there always exists such an embedding,
homotopic to a constant map. The situation changes if we wish to prescribe the
topological properties of the embedding, e.x. its action in
homology/homotopy groups. For example, if we demand that
$\varp_*:H_1(M)\to H_1(Q)$ is nonzero, then, evidently, there is a
necessary condition that the length of the shortest geodesic in the homology
class
$C\in H_1(M,\bbz)$ is not less than that of $\varp_*C$.
Using the machinery developed above, we arrive
to more obstructions of global character.

\proclaim{Theorem 3} Let $M^3$ and $Q^q$ be compact Riemannian manifolds.
Suppose $M$ is homologically atoroidal with $\pi_2(M)=0$ and $b_1(M)$,
$b_1(Q)\not=0$. Then there is a constant $C(Q)$, such that if there exists an
isometrical
immersion $f:M\to Q$ inducing a nontirivial map in the first real homology,
then
$$
R(M)\cdot\Vol(M)\ge C(Q).
$$
\endproclaim

\demo{Proof} Let $\om_1,\ldots,\om_r$ be a basis of harmonic $1$-forms
representing integer classes in $H^1(Q,\bbr)$.
Put $C^{-1}(Q)=\max\limits_{i}||\om_i||_{L^\infty}$.
If the action of $f_*:H_1(M,\bbr)\to H_1(Q,\bbr)$ is nontrivial,
then
$[f^*\om_i]\not=0$ for some $i$.
Since $[\om_i]\in H_1(Q,\bbz)$, also $[f^*\om_i]\in H_1(M,\bbz)$.
Applying the proposition 3, we get
$$
||f^*\om_i||_{L^1}\ge\frac{4\pi}{R(M)}.
$$
But
$||f^*\om_i||_{L^1}=\int_M||f^*\om_i||d\Vol\le||\om_i||_{L^\infty}\cdot\Vol(M)$,
so $R(M)\cdot\Vol(M)\ge C(Q)$. \hfill Q.E.D.

\centerline{References}

\item{[1]} T. Aubin, {\sl \'Equations diff\'erentielles non lin\'eares
et probl\'eme de Yamabe concernant la courboure scalaire}, J. Math. Pures Appl.
{\bf 55}(1976), 269--296.

\item{[2]} M. Gromov, {\sl Partial Differential Relation}, Springer-Verlag,
1986.

\item{[3]} M. Gromov, H.B. Lawson, {\sl Positive scalar curvature and the Dirac
operator on complete Riemannian manifolds}, Publ. Math. IHES, {\bf 58}(1983),
p. 295--408.

\item{[4]} D. Johnson, J. Millson, {\sl Deformation spaces associated to
compact
hyperbolic manifolds}, in: Discrete groups in Geometry and Analysis (Papers in
honor of G.D. Mostow), Progress in Math,, {\bf 67}, Birkhauser, (1987), P.
48--106.

\item{[5]} P.Li, S.-T.Yau, {\sl A new conformal invariant and its applications
to the Willmore conjecture and the first eigenvalue of compact surfaces},Inv.
Math, {\bf 69}(1982),269--291.

\item{[6]} J.M. Lee, T.H. Parker, {\sl The Yamabe problem},
Bull. AMS. {\bf 17}(1987), 37--91.

\item{[7]} A. Reznikov, {\sl Harmonic maps, hyperbolic cohomology and higher
Milnor inequalities}, to appear in Topology.

\item{[8]} A. Reznikov, {\sl Quadratic equations in groups}, preprint, 1993.

\item{[9]} A. Reznikov, {\sl Hyperbolization of homology spheres}, preprint,
1993.

\item{[10]} J. Sacks, K. Uhlenbeck, {\sl Minimal immersions of closed
Riemann surfaces}, Trans. AMS {\bf 271}(1982), 639-652.

\item{[11]} R. Schoen, {\sl Recent progress  in geometric partial differential
equations},
Proceedings of the ICM-86, AMS, 1987.

\item{[12]} R. Schoen, {\sl Conformal deformation of a Riemannian metric to
constant scalar curvature}, J. Diff. Geom., {\bf 20}(1984), 479--495.

\item{[13]} R. Schoen, S.T. Yau, {\sl Existence of incompressible minimal
surfaces
and the topology of three dimensional manifolds of positive scalar curvature},
Ann. Math. {\bf 110}(1979), 127--142.

\item{[14]} N. Trudinger, {\sl Remarks concerning the conformal deformation of
Riemannian
structures on compact manifolds}, Ann. Scuola Norm Sap. Pisa {\bf 22}(1968),
265-274.

\item{[15]} W. Thurston, {\sl A norm for the homology of 3-manifolds}, Mem.AMS,
{\bf 59},(1986), 98--130.

\item{[16]} H. Yamabe, {\sl On a deformation of Riemannian structure on compact
manifolds},
Osaka Math. J., {\bf 12}(1960), 21--37.

\item{[17]} M. Ville, {\sl Harmonic maps, second homology classes of smooth
manifolds, and bounded cohomology}, preprint (1992).

\end